\documentclass{pspum-l}

\newtheorem{theorem}{Theorem}[section]

\newtheorem{con}[theorem]{Conjecture}
\newtheorem{prob}[theorem]{Problem}
\newtheorem{proposition}[theorem]{Proposition}

\theoremstyle{definition}

\theoremstyle{remark}

\numberwithin{equation}{section}

\newcommand{\PP}{\mathbb P}
\newcommand{\NN}{\mathbb N}
\newcommand{\ZZ}{\mathbb Z}
\newcommand{\RR}{\mathbb R}

\begin{document}

\title{Some open problems on Coxeter groups and unimodality}


\author{Francesco Brenti}
\address{Dipartimento di Matematica \\
Universit\'{a} di Roma ``Tor Vergata''\\
Via della Ricerca Scientifica, 1 \\
00133 Roma, Italy \\}
\email{brenti@mat.uniroma2.it}

\subjclass[2010]{Primary 05A20, 20F55; Secondary 05E10, 05E16}

\date{}

\begin{abstract}
In this paper I present some open problems on Coxeter groups and unimodality, together with the main partial results, and computational evidence, that are 
known about them.
\end{abstract}

\maketitle

\section{Coxeter groups}

Recall that a {\em Coxeter system} is a pair $(W,S)$ where $ S := \{ s_1, \ldots , s_n  \}$ is a finite set and $W$ is a group having $S$ as generating set and relations of the form 
\[
(s_i s_j)^{m_{i,j}}=e
\]
for all $i,j \in [n]$, where $m_{i,j} \in \PP \cup \{ + \infty \}$ 
(where $\PP := \{ 1, 2, 3, \ldots \}$),
$m_{i,i}=1$, $m_{i,j}=m_{j,i} \geq 2$ if $i \neq j$, and there is no relation if
$m_{i,j}= +\infty$. Given $w \in W$ we let
\[
\ell(w) := min \{ k \in \NN : \mbox{ there are } s_{i_1}, \ldots , s_{i_k} \in S 
\mbox{ such that }  w=s_{i_1} \cdots s_{i_k} \}
\]
(where $\NN := \PP \cup \{ 0 \}$) and call this the {\em length} of $w$. The {\em right} (resp. {\em left}) {\em descent set} of $w$ is then 
\[
D_R(w) := \{ s \in S : \ell(w)>\ell(ws) \}
\]
and
\[
D_L(w) := \{ s \in S : \ell(w)>\ell(sw) \}.
\]
We let 
\[
T := \{ wsw^{-1} : s \in S, w \in W \}
\]
be the set of {\em reflections } of $W$, and
\[
N_L(w) := \{ t \in T : \ell(w)>\ell(tw) \}.
\]
Recall that the {\em Bruhat graph} of $W$ 
is the directed graph $B(W)$ having $W$ as set of vertices and where, for all $u,v \in W$, $u \rightarrow v$ if and only if $v=u \, t$ for some $t \in T$ and $\ell(v)>\ell(u)$,
in this case we also write $u \stackrel{t}{\rightarrow} v$.
The {\em Bruhat order} on $W$ is the partial order on $W$, denoted by $\leq$,
that is the transitive closure of the Bruhat graph. For $u,v \in W$ we let
$[u,v]:= \{ z \in W : u \leq z \leq v \}$ and $\ell(u,v):= \ell(v)-\ell(u)$.

It is well known (see, e.g., \cite[Prop. 5.12]{Hum} or \cite[Cor. 7.1.4]{BB}) 
that if $(W,S)$ is a Coxeter system then $\sum_{u \in W} x^{\ell(u)}$ is a rational generating function.
The following problem was proposed by Stembridge in \cite{Ste97}.
\begin{prob}
Let $(W,S)$ be a Coxeter system. Is it true that then
\[
\sum_{t \in T} x^{\ell(t)}
\]
is a rational generating function?
\end{prob}
\noindent 
No partial results on this problem seem to be known.
\medskip

The {\em Kazhdan-Lusztig polynomials} of $W$ are a family of polynomials 
$\{ P_{u,v}(q)\}_{u,v \in W}$  $\subseteq \ZZ[q]$ in one variable $q$, indexed by
pairs of elements of $W$, that were first defined by Kazhdan and Lusztig in \cite{KL}.
These polynomials have found important applications in various areas of mathematics (see, e.g., \cite{KaLu80}, \cite{BeiBer}, \cite{BryKas}, and also \cite{BilLak}, \cite{BB}, \cite{HumBGG}, and the references cited there) and have been generalized in
various ways (see, e.g., \cite{superKL}, \cite[\S 6]{StaJAMS}, \cite{BreKLS}, \cite{ParabolicKL},
\cite{LusVog}, \cite{MatroidKL}, \cite{DiamondKL}, \cite{KLTL}, to cite just a few).
We find it most convenient here to define the parabolic Kazhdan-Lusztig polynomials.
Let $J \subseteq S$. Recall that the (left) quotient of $W$ corresponding to
$J$ is 
\[
W^J := \{ u \in W : D_L(u) \subseteq S \setminus J \}.
\]
The polynomials can be defined in two ``Theorem-Definitions'' (\cite{ParabolicKL}).
\begin{theorem}
	\label{DefR}
	Let $(W,S)$ be a Coxeter system, $J \subseteq S$, and $x \in \{ -1,q \}$.
	Then
	there is a unique family of polynomials $\{ R_{u,v}^{J,x}
	\} _{u,v \in W^{J}} \subseteq {\mathbb Z}[q]$ such that, for all $u,v \in W^{J}$:
	\begin{description}
		\item[i)] $R^{J,x}_{u,v}=0$ if $u \not \leq v$;
		\item[ii)] $R^{J,x}_{u,u}=1$;
		\item[iii)] if $u<v$ and $s \in D_R(v)$ then
		\[ R_{u,v}^{J,x} (q)= \left\{ \begin{array}{ll}
		R_{us,vs}^{J,x}(q), & \mbox{if $us<u$,} \\
		(q-1)R_{u,vs}^{J,x}(q)+qR_{us,vs}^{J,x}(q), & \mbox{if $u<us \in W^{J}$,} \\
		(q-1-x)R_{u,vs}^{J,x}(q), & \mbox{if $u<us \not \in
			W^{J}$.}
		\end{array} \right. \]
	\end{description}
\end{theorem}

\begin{theorem}
	\label{DefKL}
	Let $(W,S)$ be a Coxeter system,  $J \subseteq S$, and $x \in \{ -1,q \}$.
	Then 
	there is a unique family of polynomials $\{ P^{J,x}_{u,v} \}
	_{u,v \in  W^{J}}  \subseteq {\mathbb Z} [q]$ such that, for all
	$u,v \in W^{J}$:
	\begin{description}
		\item[i)] $P^{J,x}_{u,v}=0$ if $u \not \leq v$;
		\item[ii)] $P^{J,x}_{u,u}=1$;
		\item[iii)] deg$(P^{J,x}_{u,v}) <  \frac{1}{2}
		\ell (u,v) $ if $u < v$;
		\item[iv)]
		\[ q^{\ell (u,v)} \, P^{J,x}_{u,v} \left( \frac{1}{q} \right)
		= \sum _{a \in [u,v]^{J}}   R^{J,x}_{u,a}(q) \,
		P^{J,x}_{a,v}(q) \]
		if $u \leq v$, where
		\[
		[u,v]^J := \{ a \in W^J : u \leq a \leq v \}.
		\]
	\end{description}
\end{theorem}
\noindent 
The polynomials $\{ R_{u,v}^{J,x} \} _{u,v \in W^{J}} $ (resp., $\{ P_{u,v}^{J,x}
\} _{u,v \in W^{J}}$) whose existence the previous results show are called  the {\em
	parabolic $R$-polynomials} (resp., {\em parabolic  Kazhdan-Lusztig polynomials})
of $W^{J}$ of type $x$. The polynomials 
$R_{u,v}:=R_{u,v}^{\emptyset ,-1}$ ($=R_{u,v}^{\emptyset ,q}$) and 
$P_{u,v}:=P_{u,v}^{\emptyset ,-1}$ ($=P_{u,v}^{\emptyset ,q}$) for $u,v \in W$,
are called the {\em $R$-polynomials} (resp., {\em Kazhdan-Lusztig polynomials})
of $W$.

The most outstanding open problem about Kazhdan-Lusztig polynomials, particularly 
from a combinatorial point of view, is the so-called ``Combinatorial Invariance
Conjecture''.
\begin{con}
	Let $(W_1,S_1)$, $(W_2,S_2)$ be two Coxeter systems, and $u,v \in W_1$, $w,z \in W_2$ be such that $[u,v] \simeq [w,z]$ (isomorphic as posets). Then
	\begin{equation}
	\label{CIC}
	P_{u,v}(q)=P_{w,z}(q).
	\end{equation}
\end{con}

This conjecture was made by Lusztig \cite{BjoPC} (see also \cite[Rem. 7.31]{Dye87}).
It may be noted that Kazhdan never conjectured the CIC. However, during a long
mathematical conversation in Jerusalem \cite{KazPC01} he told me: ``I believe the
[combinatorial invariance] conjecture to be true''. 
The conjecture is known to be true if $\ell(u,v) \leq 4$ (\cite[7.31]{Dye87},  \cite[Chap. 5, Exercises 7 and 8]{BB}), if $\ell(u,v)\leq 8$ and $W_1$ and $W_2$ are both of type
$A$ (\cite{Inc07}), if $\ell(u,v)\leq 6$ and $W_1$ and $W_2$ are both of type
$B$ or $D$ (\cite{Inc07}), 
if $W_1$ and $W_2$ are both of type $\tilde{A}_2$ (\cite{BLP}), 
if $[u,v]$ is a lattice (\cite[7.23]{Dye87}, \cite[Thm. 6.3]{BreInv}), 
if one adds the hypothesis that ``$P_{u,v}(q)=1$'' (\cite[Thm. C]{Car94}, \cite[Prop. 3.3]{DyeBru}), 
if $u v^{-1} \in T$ and $W_1$ and $W_2$ are both of type $A$ (\cite[Cor. 4.7]{BreAdv},
\cite[Prop. 3.3]{DyeBru}), 
if $x=u=e$ (\cite{BCM}), and for the coefficient of $q$ if $W_1$ and $W_2$ are both simply laced (\cite{Patimo}). Note that an equivalent conjecture is
obtained if one replaces (\ref{CIC}) with ``$R_{x,y}(q)=R_{u,v}(q)$''. It would be 
desirable to have computational evidence in favor of the CIC, particularly for the
exceptional finite Coxeter groups. I have been asked more than once if I believe the CIC. If I had to bet my money, I would say that it's true.
The following result, proved in \cite[Prop. 3.3]{DyeBru}, has a flavor similar to the CIC, and plays a role in some of the proofs of the CIC in special cases.
For $u,v \in W$ let $B(u,v)$ be the directed graph induced on $[u,v]$ by $B(W)$
(so, the vertex set of $B(u,v)$ is $[u,v]$ and, if $x,y \in [u,v]$, there is
a directed edge from $x$ to $y$ in $B(u,v)$ if and only if $x \rightarrow y$).
\begin{theorem}
	Let $(W_1,S_1)$, $(W_2,S_2)$ be two Coxeter systems, and $u,v \in W_1$, $w,z \in W_2$ be such that $[u,v] \simeq [w,z]$ (isomorphic as posets). Then
	$B(u,v) \simeq B(w,z)$ (isomorphic as directed graphs).
\end{theorem}

The most natural generalization of the CIC to the parabolic polynomials is that if 
$(W_1,S_1)$ and $(W_2,S_2)$ are two Coxeter systems, $J_1 \subseteq S_1$, $J_2 \subseteq S_2$, and $u,v \in W_1^{J_1}$, $w,z \in W_2^{J_2}$ are such that $[u,v]^{J_1} \simeq [w,z]^{J_2}$ (as posets) then $P_{u,v}^{J_1,q}(q)=P_{w,z}^{J_2,q}(q)$ 
(equivalently, $P_{u,v}^{J_1,-1}(q)=P_{w,z}^{J_2,-1}(q)$, $R_{u,v}^{J_1,q}(q)=R_{w,z}^{J_2,q}(q)$, $R_{u,v}^{J_1,-1}(q)=R_{w,z}^{J_2,-1}(q)$).
Although this statement holds for some particularly nice and interesting quotients
(\cite{BrePac}, \cite{BreTAMS}) it's false in general (\cite{BMS}). For example,
if $W_1$ and $W_2$ are both of type $B_5$, $J_1=J_2=S \setminus \{ s_3 \}$ (numbering as in \cite[Appendix A1]{BB}), $u=[4,1,5,2,3]$, $v=[5,-4,1,2,3]$,
 $w=[1,4,2,5,3]$, and $z=[-4,1,5,2,3]$ then $[u,v]^{J_1} \simeq [w,z]^{J_2}$ while
 $P_{u,v}^{J_1,q}(q)=q \neq 0 = P_{w,z}^{J_2,q}(q)$.
 
 A deeper generalization of the CIC to the parabolic setting has been proposed
 by Marietti in \cite{MaTAMS}.
 \begin{con}
 	\label{Par_CIC}
 	Let $(W_1,S_1)$, $(W_2,S_2)$ be two Coxeter systems, $J_1 \subseteq S_1$, $J_2 \subseteq S_2$, $u,v \in W_1^{J_1}$, $w,z \in W_2^{J_2}$ and $f: [u,v] \rightarrow [w,z]$ be a poset isomorphism such that $f([u,v]^{J_1})=[w,z]^{J_2}$. Then 
 	\begin{equation}
 	\label{ParCIC}
 	P_{u,v}^{J_1,q}(q)=P_{w,z}^{J_2,q}(q).
 	\end{equation}
 \end{con}

Again, an equivalent statement is obtained by substituting (\ref{ParCIC}) with
``$P_{u,v}^{J_1,-1}(q)$ $=P_{w,z}^{J_2,-1}(q)$'' or the analogous equalities for the 
parabolic $R$-polynomials. It is clear that Conjecture \ref{Par_CIC} reduces to the
CIC if $J_1=J_2=\emptyset$, and that it holds if the naive parabolic CIC does. 
In addition, Conjecture \ref{Par_CIC} is known to hold if $u=w=e$ (\cite{MaAM}).

Just as finding a combinatorial interpretation is more satisfactory than 
proving nonnegativity so finding an explicit algorithm for computing
the Kazhdan-Lusztig (or $R$) polynomial of a pair of elements $u,v \in W$
starting from the Bruhat interval $[u,v]$ as an abstract poset 
would be more satisfactory than simply proving the CIC.
Recently, using techniques from deep learning (\cite{Nat}), a candidate such algorithm has been proposed, in the case that $W$ is a Weyl
group of type $A$, in \cite{Wil+}. 
Let $u,v \in W$, $u \leq v$. For any $z \in [u,v]$ $z < v$, say that $[u,z]$ is
{\em diamond complete} in $[u,v]$ if whenever $a,b,c,d \in [u,v]$ are such that
$a \rightarrow b$, $b \rightarrow d$, $a \rightarrow c$, $c \rightarrow d$, 
and $a,b,c \in [u,z]$ then $d \in [u,z]$. Let $[u,z]$ be such a diamond complete
subinterval. For any $w \in [u,z]$ let $U(w) := \{ x \in [u,v] \setminus [u,z] :
w \rightarrow x \}$, and consider the partial order induced on $U(w)$ by the
Bruhat order. Say that $U(w)$ {\em spans a hypercube cluster} if whenever
$A \subseteq U(w)$ is an antichain then there is a unique
embedding of directed graphs $\theta : {\mathcal P}(A) \rightarrow B(u,v)$ such
that $\theta( \emptyset)=w$ and $\theta ( \{ a \})= a$ for all $a \in A$
(where ${\mathcal P}(A)$ is the directed graph having the subsets of $A$ as vertices and where $B \rightarrow C$ if and only if $B \subset C$ and $|C \setminus B|=1$).
Finally, say that $[u,z]$ is a {\em hypercube decomposition} 
if ($[u,z]$ is diamond complete and) $U(w)$ spans a hypercube cluster
for any $w \in [u,z]$. In \cite{Wil+} a procedure is given 
to compute, starting
from a hypercube decomposition $[u,z]$ of a Bruhat interval $[u,v]$, and all the Kazhdan-Lusztig polynomials $P_{x,y}(q)$ for $x,y \in [u,v]$ such that 
$\ell(x,y)<\ell(u,v)$, a polynomial $\widetilde{P}_{u,v,z}(q)$, and the following conjecture
is made (\cite[Conj. 3.8]{Wil+}).
\begin{con}
\label{Wilcon} 
	Let $W$ be a Weyl group of type $A$ and $u,v \in W$, $u \leq v$. Then
\begin{equation}
\label{Wilconj}
P_{u,v}(q) = \widetilde{P}_{u,v,z}(q)
\end{equation}
for any hypercube decomposition $[u,z]$ of $[u,v]$.
\end{con}
\noindent 
It is shown in \cite[Thm. 3.7]{Wil+} that if $u,v \in S_n$ then
$\{ w \in [u,v] : w^{-1}(1)=u^{-1}(1) \}$ is a hypercube decomposition
of $[u,v]$ and (\ref{Wilconj}) holds for this choice. 
Conjecture \ref{Wilcon} is true for all intervals in $S_n$ if $n \leq 7$, 
and for millions of intervals in $S_8$ and $S_9$. Note that a  
Bruhat interval in a general Coxeter group may not have a hypercube decomposition
(for example, a $5$-crown does not have a hypercube decomposition).
\medskip

Although the nonnegativity conjecture for the Kazhdan-Lusztig polynomials (\cite{KL})
has been proved (\cite{EW}) there is (at least) one other nonnegativity conjecture
that is still open. Recall (see \cite[\S 2]{Dye93}) that a {\em reflection ordering} on $(W,S)$ is a total order $\preceq$ on $T$ such that if $W'$
is a dihedral reflection subgroup of $W$ (so $W':= <J>$ for some $J \subseteq T$ 
and $S':=|\{ t \in T : N_L(t) \cap W' = \{ t \} \}|=2$) then either
\[
a \preceq aba \preceq ababa \preceq \cdots \preceq babab \preceq bab \preceq b
\]
or 
\[
a \succeq aba \succeq ababa \succeq \cdots \succeq babab \succeq bab \succeq b
\]
where $\{ a,b \}:= S'$. Equivalently, using the canonical bijection between $T$ and $\Phi^+$ (see, e.g., \cite[Prop. 4.4.5]{BB}), a reflection ordering is a total order on $\Phi^+$ such that if $\alpha,\beta \in \Phi^+$ and $\lambda, \mu \in \RR_{>0}$ are such that $\lambda \alpha + \mu \beta \in \Phi^+$
then either $\alpha \preceq \lambda \alpha + \mu \beta \preceq \beta$ or 
$\alpha \succeq \lambda \alpha + \mu \beta \succeq \beta$. 
Reflection orderings always exist, and 
in fact there are many (we refer the reader to \cite[\S 5.2]{BB}, and \cite{Dye93} for
further information about reflection orderings). If $W$ is of type $A$ then the lexicographic order $(1,2) \prec (1,3) \prec \cdots \prec (1,n) \prec (2,3) \prec \cdots \prec (n-1,n)$ is a reflection ordering.

Let $\prec$ be a reflection ordering of $(W,S)$ and $\Gamma = x_0 \stackrel{t_1}{\rightarrow} x_1 \stackrel{t_2}{\rightarrow} \cdots \stackrel{t_k}{\rightarrow} x_k$ be a directed path in the Bruhat graph $B(W)$.
Let $a$ and $b$ be two noncommuting variables. We associate to $\Gamma$ a monomial $m_{\prec}(\Gamma) := y_1 \cdots y_{k-1}$ in $a$ and $b$ by letting
\[ 
y_j := \left\{ \begin{array}{ll}
a, & \mbox{if $t_j \prec t_{j+1}$,} \\
b, & \mbox{if $t_j \succ t_{j+1}$,}
\end{array} \right. 
\]
for $j=1, \ldots, k-1$, and $m_{\prec}(\Gamma) := 1$ if $k=1$, and let
\[
\widetilde{\psi}_{u,v} := \sum_{\Gamma} m_\prec(\Gamma)
\]
where the sum is over all directed paths $\Gamma$ from $u$ to $v$ in $B(W)$. 
It can be shown
(see \cite[Prop. 1.5]{BilBre}, and also \cite[Prop. 4.4]{BreLMS}) that, although
$m_\prec(\Gamma)$ depends on the reflection ordering used to define it,
$ \widetilde{\psi}_{u,v}$ does not. Also, it is known (see, e.g., the proof of
Theorem 4 in \cite{BayKla}) that there is a polynomial 
$ \widetilde{\Phi}_{u,v} \in \ZZ \langle c,d \rangle$ in two non-commuting variables $c$ and $d$ such that
\[
 \widetilde{\psi}_{u,v}(a,b) =  \widetilde{\Phi}_{u,v}(a+b,ab+ba).
\]
The polynomial $ \widetilde{\Phi}_{u,v}(c,d)$ is called the {\em complete cd-index}
of $u,v$. The reason for this terminology lies in the fact that the homogeneous part
of highest possible degree in $ \widetilde{\Phi}_{u,v}(c,d)$ (where $deg(a)=deg(b)=deg(c)=1$ and $deg(d)=2$) namely $\ell(u,v)-1$ is (by \cite[Thm. 3.14.2]{StaEC1} and the 
fact that the assignment $\lambda(w,z):= w^{-1} z$ for $u \leq w \lhd z \leq v$ is 
an $EL$-labeling of $[u,v]$, \cite{Dye93}) the $cd$-index of $[u,v]$ as an Eulerian poset
(see, e.g., \cite[\S 3.17]{StaEC1}). It is known (\cite[Thm. 4.1]{BilBre}) that the
coefficients of the Kazhdan-Lusztig polynomial $P_{u,v}(q)$ are given by explicit linear combinations of the coefficients of $ \widetilde{\Phi}_{u,v}(c,d)$.
The following conjecture appears in \cite[Conj. 6.1]{BilBre}.
\begin{con}
	\label{compcd}
	Let $(W,S)$ be a Coxeter system, and $u,v \in W$, $u \leq v$. Then
	\[
	 \widetilde{\Phi}_{u,v} \in \NN \langle c,d \rangle.
	\]
\end{con}
Conjecture \ref{compcd} has been verified if $W$ is of type $A$ and $\ell(u,v) \leq 7$, 
and is known to be true for the coefficients of the monomials of highest possible degree
by \cite[Thm. 1.3]{Karu} since every Bruhat interval $[u,v]$ is a Gorenstein$^*$
poset. Some further evidence is presented in \cite[\S 6]{BilBre} and \cite{Blanco}.

It is known (\cite{Polo}, see also \cite{CasRep}) that if $P(q) \in \NN[q]$ is such that $P(0)=1$ then there are $n \in \PP$ and $u,v \in S_n$ such that $P(q)=P_{u,v}(q)$. The following related problem was posed by Bj\"{o}rner (\cite{BjoPC}).
\begin{prob}
Let $(W,S)$ be a Coxeter system and $u,v \in W$, $u \leq v$. Is it true that
then there are a Coxeter system $(W',S')$ and $w \in W'$ such that 
\[
P_{u,v}(q)=P_{e,w}(q) \; ?
\]
\end{prob}
\noindent 
It is easy to see that the answer to this question is negative if one requires
that $W'=W$. No partial results on this problem seem to be known.

The proof of the celebrated nonnegativity conjecture for Kazhdan-Lusztig polynomials 
(\cite{EW}) makes the following problem even more compelling.
\begin{prob}
	\label{combintprob}
	Find a combinatorial interpretation for Kazhdan-Lusztig polynomials.
\end{prob}
\noindent
So, given a Coxeter system $(W,S)$ and $u,v \in W$, $u \leq v$, one would like to produce (in some explicit combinatorial way) a set $M(u,v)$ and function $s: M(u,v) \rightarrow \NN$ such that
\begin{equation}
\label{combint}
P_{u,v}(q) = \sum_{a \in M(u,v)} q^{s(a)}.
\end{equation}
Problem \ref{combintprob} is open, and interesting, even for $q=1$. Combinatorial
interpretations for the Kazhdan-Lusztig polynomials are known if $W$ is a Weyl group and $v$ is a Deodhar element (\cite[Thms. 2.3 and 5.12]{BilJon}), if $(W,S)$ is a universal
Coxeter system (\cite[Thm. 3.8]{DyeJA}), if $W$ is a Weyl group and $v$ is rationally smooth (\cite[Thm. 2.5]{BilPos}, \cite[Thm. A2]{KL}), if the Coxeter graph of $(W,S)$ is acyclic and $u$ and $v$ are Boolean elements (\cite[Cor. 4.3]{Mon}), if the 
Coxeter graph of $(W,S)$ is a cycle and $u$ and $v$ are Boolean elements
(\cite[Thm.4.4]{Mar}), if $W$ is a Weyl group and $u,v \in W^J$ where $(W,W_J)$ is
a Hermitian symmetric pair ($J \subseteq S$) (\cite{Boe}, see also \cite{EnrShe}),
and if $[u,v]$ is isomorphic (as a poset) to a Boolean algebra or to the lattice of faces
of an $(\ell(u,v) -1)$-dimensional cube (\cite[Cor. 6.8 and 6.9]{BreInv}). In addition to the above cases, if $W$ is of type $A$ then combinatorial interpretations are known for various special families of permutations (\cite{CasJAC}, \cite{CasRep}, \cite{CasMarDM}, \cite{Woo}), including if $u=e$, $v([3]) =[n-2,n]$, $v([n-2,n])=[3]$
and $v(4)>v(5)> \cdots >v(n-3)$ (\cite[Cor. 5.5]{CasMarDM}, \cite[Thm. 4.8]{CasJAC}).
The combinatorial interpretation given in \cite{LasCR95} (see also \cite[Chap.5, Ex.39]{BB}) in the case that $v$ is a permutation that avoids 3412 must be considered a conjecture since no complete proof of that statement is known. It should be noted that
Deodhar in \cite{DeoGD} constructs, given a Coxeter system $(W,S)$ for which the 
Kazhdan-Lusztig polynomials have nonnegative coefficients and $u,v \in W$, $u \leq v$,
a set  $M(u,v)$ and function $s: M(u,v) \rightarrow \NN$ such that (\ref{combint}) holds. However, the definition of $M(u,v)$ is recursive and for this reason this is not generally considered a ``combinatorial interpretation''. Still, the combinatorial interpretations for Deodhar elements referred to above, and a few others, have been obtained using Deodhar's general framework. For a more algebraic viewpoint on the combinatorial interpretation problem, and Deodhar's construction, see \cite{LibWil}.
A hint for a general combinatorial interpretation could come from the following result, which is proved in \cite[Thm. 5.8]{Pla}.
\begin{theorem}
	Let $(W,S)$ be a Coxeter system and $u,v,w \in W$, $u \leq v \leq w$. Then
	\[ P_{v,w}(q) \leq P_{u,w}(q) \]
	(coefficientwise).
\end{theorem}

Just as for the CIC one may consider the combinatorial interpretation problem also for the parabolic Kazhdan-Lusztig polynomials. The nonnegativity of the parabolic Kazhdan-Lusztig
polynomials of type $q$, $\{ P^{J,q}_{u,v}(q) \}_{u,v \in W^J}$, is proved in \cite[Thm. 1.1]{LibWilpre}. Combinatorial interpretations of these polynomials are known if $W$ is of type $A$ and $W^J$ is a tight quotient (\cite{BrePac}, \cite{BIM}), if $W$ is a Weyl group and $W^J$ is a quasi-minuscule quotient (\cite{BrePac}, \cite{BreTAMS}, \cite{BMS}, see also \cite{ShiZin}, \cite{LejStro}, \cite{CoxVis}), and if $u$ and $v$ are Boolean elements and $W$ 
is of type $A$ (\cite[Thm. 5.2]{MarEJC}). Also, it is known
(\cite[Thm. 5.8]{Sen}) that it is enough to find a combinatorial interpretation for the parabolic Kazhdan-Lusztig polynomials $\{ P^{J,x}_{u,v} \}_{u,v \in W^J }$ in the case that
$|J|=|S|-1$. A hint for a combinatorial interpretation could come from the following result which was conjectured in  \cite{BreStoc08} and proved in \cite[Cor. 8.4]{LibWilpre}.
\begin{theorem}
Let $(W,S)$ be a Coxeter system and $I \subseteq J \subseteq S$. Then
\[ 
P^{J,q}_{u,v}(q) \leq P^{I,q}_{u,v}(q) 
\]
(coefficientwise) for all $u,v \in W^J$.
\end{theorem}
The nonnegativity  of the parabolic Kazhdan-Lusztig polynomials of type $-1$ is open in general, but see remark (2) after Theorem 1.1 in \cite{LibWilpre}.

\section{Unimodality}

Recall that a sequence $(a_0, \ldots, a_n)$ is said to be {\em unimodal} if 
there is $0 \leq m \leq n$ such that 
$a_0 \leq \cdots \leq a_{m-1} \leq a_m \geq a_{m+1} \geq \cdots \geq a_n$, 
and is said to be
{\em symmetric} if $a_j = a_{n-j}$ for all $j=0, \ldots , n$. It is said to 
be {\em log-concave}
if $(a_j)^2 \geq a_{j-1} a_{j+1}$ for all $1 \leq j \leq n-1$ and 
{\em ultra log-concave} if the sequence 
$\{ \frac{a_j}{{n \choose j}} \}_{j=0,\ldots,n}$ is log-concave.
We say that a polynomial $\sum_{j=0}^{n} a_j t^j$ is {\em unimodal} (resp. {\em symmetric}, {\em log-concave}, {\em ultra log-concave}) if the sequence $(a_0, \ldots, a_n)$ has the corresponding property.
If $P(t)=\sum_{j=0}^{n} a_j t^j$ is symmetric then there are unique numbers $\gamma_0, \ldots , \gamma_{\lfloor n/2 \rfloor}$ such that
\[
P(t) = \sum_{k=0}^{\lfloor n/2 \rfloor} \gamma_k \, t^k (1+t)^{n-2k}.
\]
The vector $(\gamma_0, \ldots , \gamma_{\lfloor n/2 \rfloor})$ is called the {\em $\gamma$-vector} of $P(t)$, and the polynomial $P(t)$ is said to be 
{\em $\gamma$-nonnegative} if $\gamma_k \geq 0$ for all $0 \leq k \leq \lfloor n/2 \rfloor$. It is clear that a $\gamma$-nonnegative polynomial is unimodal, and it
is not hard to see (see, e.g., \cite[Rem. 1.3.1]{BraUni}) that if $P(t)$ is a symmetric polynomial with nonnegative coefficients and only real roots then 
$P(t)$ is $\gamma$-nonnegative.  Unimodal, log-concave, and $\gamma$-nonnegative polynomials appear often in combinatorics, geometry, and algebra (see, e.g., 
\cite{StaUni}, \cite{BreUni}, \cite{BraUni}, \cite{AthUni}, and the references cited there).
A real (finite or infinite) matrix $A:= ( A_{i,j} )_{i,j \in \NN}$
is said to be {\em totally positive} (or, {\em TP}, for short, sometimes more
appropriately called {\em totally nonnegative}, \cite{FJ})
if all the minors of $A$ have nonnegative determinant. 
There is a deep relationship between polynomials with only real roots and totally
positive matrices. The following result was first proved in 
\cite{Edrei} (see also \cite[Thm. 4.5]{Pin}).
\begin{theorem}
\label{Edr}
Let $P(t)=\sum_{j=0}^{n} a_j t^j$ be a polynomial with nonnegative
coefficients. Then $P(t)$ has only real roots if and only if the matrix
\[
(a_{j-i})_{i,j \in \NN}
\]
is totally positive (where $a_k :=0$ if $k<0$ or $k>n$).
\end{theorem}
\noindent
Totally positive matrices
arise in a number of areas in science, including mathematics, statistics, probability, economics,
mechanics, and computer science (see, e.g., \cite{Karl}, the Foreword of \cite{Pin}, \cite[Sec. 0.2]{FJ} and the 
references cited there).

\medskip

For $\sigma \in S_n$ let
\[
L(\sigma) := | \{ (i,j) \in [n]^2 : i<j, \, \sigma(i)>\sigma(j), \; i \not \equiv j \pmod{2} \} |.
\]
The statistic $L$ is known as the {\em odd length} or {\em odd inversion number}
of $\sigma$. The statistic was introduced in \cite{KloVol} in relation to formed spaces
and has been further studied in \cite{BreCav} and \cite{BreCa2}.
The following conjecture appears in \cite[Conj. 6.2]{BreCa2}.
\begin{con}
	Let $n \in \PP$, $n \geq 5$. Then 
	\[
L_n(x) := \sum_{\sigma \in S_n} x^{L(\sigma)}
	\]
	is symmetric and unimodal.
\end{con}
The conjecture has been verified for $n \leq 11$. The symmetry statement is clear since
$L(\sigma w_0)=L(w_0)-L(\sigma)$ for all $\sigma \in S_n$ where
$w_0$ is the longest permutation $w_0=n \, n-1 \cdots 3 \, 2 \, 1$. Note that, in
general, $L_n(x)$ is not log-concave and not $\gamma$-nonnegative.

\medskip

Let $(W,S)$ be a Coxeter system. The following result follows easily from the 
case $J=\emptyset$ of Theorem \ref{DefR}.
\begin{proposition}
	Let $(W,S)$ be a Coxeter system. Then there exists a unique family of polynomials
	$\{ \widetilde{R}_{u,v}(t) \}_{u,v \in W} \subseteq \NN[t]$ such that
	\[
	R_{u,v}(q)= q^{\ell(u,v)/2}  \widetilde{R}_{u,v}(q^{1/2}-q^{-1/2}).
	\]
\end{proposition}
\noindent 
So, knowledge of the $ \widetilde{R}$-polynomials is equivalent to knowledge of
the $R$-polynomials. Combinatorial interpretations of the polynomials 
$\{  \widetilde{R}_{u,v}(t) \}_{u,v \in W}$ have been given in \cite{DeoInv} and 
\cite{Dye93}, see \cite[Thms 5.3.4 and 5.3.7]{BB}.
It is easy to see that, if $u \leq v$, then $deg( \widetilde{R}_{u,v})=\ell(u,v)$ and the powers 
appearing in $ \widetilde{R}_{u,v}(t)$ are all of the same parity. More precisely,
there is a polynomial $Q_{u,v}(t) \in \NN[t]$ such that $ \widetilde{R}_{u,v}(t)=
Q_{u,v}(t^2)$ if $\ell(u,v) \equiv 0 \pmod{2}$ and  $ \widetilde{R}_{u,v}(t)=
t \, Q_{u,v}(t^2)$ if $\ell(u,v) \equiv 1 \pmod{2}$. It is of interest to 
characterize the $R$-polynomials. In this respect I propose the following conjecture
which is a generalization of \cite[Conj. 7.1]{BreDM}.
\begin{con}
	Let $(W,S)$ be a finite Coxeter system and $u,v \in W$, $u \leq v$. Then the polynomial
	$Q_{u,v}(t)$ is log-concave.
\end{con}
The conjecture has been verified if $W$ is of type $F_4$, or $H_3$, or 
$A_n$, or $B_n$, or $D_n$ and $n \leq 5$, and for dihedral groups.
Note that the polynomial $Q_{u,v}(t)$ is not, in general, ultra log-concave.
For example, if $W$ is of type $A_5$, $u=213465$, and $v=563412$ then
$ \widetilde{R}_{u,v}=t^2+4t^4+6t^6+5t^8+t^{10}$ and neither the sequence $(0,1,4,6,5,1)$
nor the sequence $(1,4,6,5,1)$ are ultra log-concave.
It is conceivable that the polynomials $Q_{u,v}(t)$ are log-concave for
any Coxeter system.
\medskip

Despite the settling (\cite{Huh}) of the celebrated log-concavity 
conjecture for chromatic polynomials (\cite{Read}) there is (at least) one
other unimodality statement about them that is still open.
Let $G=(V,E)$ be a (simple, loopless) graph on $p$ vertices and $\chi(G;x)$ be 
its chromatic polynomial. Write
\[
\chi(G;x) = \sum_{i=0}^{p} (-1)^{p-i} c_i \langle x \rangle_i
\]
where $\langle x \rangle_i := x(x+1) \cdots (x+i-1)$ for $i \geq 1$ and $\langle x \rangle_0 :=1$. The {\em $\tau$-polynomial} of $G$ is
\[
\tau(G;x):= \sum_{i=0}^{p} c_i x^i.
\]
The $\tau$-polynomial was first introduced and studied in \cite{BreChro} where the following combinatorial interpretation of it is given. Let $\pi \in \Pi(V)$ where
$\Pi(V)$ is the set of all set partitions of
$V$, say $\pi = \{  B_1, \ldots , B_k \}$. For $B \subseteq V$ let $G[B]$ be the 
subgraph induced by $B$ (so $G[B]=(B, \{ \{i,j\} \in E: i,j \in B \})$) and  
$G[\pi] := \biguplus_{i=1}^{k} G[B_i]$. Finally, for a graph $H$, let $a(H)$ be the
number of acyclic orientations of $H$.
\begin{theorem}
	Let $G$ be a graph. Then
\[
\tau(G;x) = \sum_{\pi \in \Pi(V)} a(G[\pi]) \, x^{|\pi|}.
\]
\end{theorem}
\noindent
The following problem was  
raised in \cite[Prob. 7.1]{BreChro}.
\begin{prob}
	Does $\tau(G;x)$ have only real roots for all graphs $G$?
\end{prob}
\noindent
The answer is yes for all connected graphs on $\leq 8$ vertices, and
if the chromatic polynomial of $G$ has only real roots (\cite[Thm. 6.1]{BreChro}).

\medskip

Let $P$ be a convex polytope of dimension $d$. 
Recall (see, e.g., \cite[Chap. III, \S 1]{StaCCA})
that $P$ is {\em simplicial} if every proper face of it is a simplex. The
{\em boundary complex} $\Delta(P)$ of $P$ is then the simplicial complex of all the
proper faces of $P$. Being a simplicial complex, $\Delta(P)$ has an $h$-vector $(h_0(\Delta(P)), \ldots , h_d(\Delta(P)))$ (we refer the reader to, e.g., \cite[Chap. II, \S 2, p. 58]{StaCCA}
for the definition of the $h$-vector of a simplicial complex).
The following famous result is well known (see, e.g., \cite[Chap. III, Thm. 1.1]{StaCCA}).
\begin{theorem}
	Let $P$ be a simplicial convex polytope. Then the $h$-vector of $\Delta(P)$ 
	is symmetric and unimodal.
\end{theorem}

For the barycentric subdivision of a simplicial convex polytope more can be said 
(\cite[Cor. 3]{BreWel}).
	
\begin{theorem}
	Let $P$ be a simplicial convex polytope. Then the (generating polynomial of the) $h$-vector of the barycentric subdivision of $P$ has only real roots.
\end{theorem}

The following problem is natural, and has been circulated informally by the author
since 2004.
\begin{prob}
	\label{SCPHTP}
	Let $\sum_{i=0}^{d} h_i t^i \in \NN[t]$ be a symmetric, monic polynomial with only real roots. Is there a simplicial convex polytope whose $h$-vector is $(h_0, \ldots , h_d)$?
\end{prob}
\noindent 
The answer is yes if $d \leq 9$ and $h_i \leq 100$ for all $0 \leq i \leq d$. 
Note that there are numerical characterizations both of
$h$-vectors of simplicial convex polytopes (\cite[Chap. III, Thm. 1.1]{StaCCA}) as well
as of monic polynomials with only real roots and nonnegative coefficients (see
Theorem \ref{Edr}) thus Problem \ref{SCPHTP} is really asking whether one such 
set of inequalities implies the other one.
The following related problem appears in \cite[Question 4.4]{KubWel}.
\begin{prob}
	\label{KW}
	Let $\sum_{i=0}^{d} h_i t^i \in \NN[t]$ be a polynomial with only real roots
	such that $h_0=1$ and $h_0 < h_1 < \cdots < h_k$ for some $0 \leq k \leq d$. Is there a simplicial complex whose $f$-vector is $(h_0, h_1-h_0, \ldots , h_k-h_{k-1})$?
\end{prob}
\noindent 
Note that, by \cite[Chap. III, Thm. 1.1]{StaCCA}, a positive answer to Problem \ref{KW} implies a positive answer to Problem \ref{SCPHTP}. Related results and problems also appear in \cite{BelSka}. 
\medskip

Let $(W,S)$ be a Coxeter system. It is of interest, and difficult, to obtain 
properties of the
rank generating function of Bruhat intervals $[u,v]$, $u,v \in W$ (see, e.g., \cite{BjoEke}).
In this respect, I feel that the following holds.
\begin{con}
	Let $W$ be a Weyl group, and $u,v \in W$. Then $[u,v]$ is rank log-concave.
\end{con}
\noindent 
The conjecture has been verified if $W$ is of type $A_n$ and $n \leq 5$, 
or $D_n$ and $n \leq 5$, or $B_n$ and $n \leq 4$, or $B_5$ and $\ell(u,v) \geq 20$, 
or $F_4$, and for the dihedral groups. Note that the corresponding statement 
does not hold for finite Coxeter systems. For example, if $(W,S)$ if of type $H_3$,
$u=s_3$, and $v=s_1s_2s_3s_2s_1s_2s_1s_3$ (where $S=\{ s_1,s_2,s_3 \}$, 
$m(s_1,s_2)=5$, and $m(s_2,s_3)=3$) then the rank generating function of $[u,v]$
is $1+3t+5t^2+7t^3+10t^4+10t^5+5t^6+t^7$.

\medskip
Many matrices arising in combinatorics
are known to be TP (see, e.g., \cite{BreMem}, \cite{BreCombTP}, and
\cite{BreTPComb}). It is therefore surprising that for a fundamental and
old combinatorial matrix such as the one consisting of the Eulerian numbers this property has not yet been settled.  
For $n \in \PP$ and $k \in \NN$ let $A(n,k)$ be the corresponding
{\em Eulerian number} (so, $A(n,k)$ is the number of permutations in $S_n$ that
have $k$ descents). The following conjecture was first put forward in \cite[Conj. 6.10]{BreCombTP}.
\begin{con} \label{Euler}
The matrix
\[
A:=( A(n+1,k))_{n,k \in \NN}
\]
is totally positive.
\end{con}
\noindent 
It has been checked that $( A(n+1,k))_{0 \leq n,k \leq 44}$ is TP.
A more general conjecture, which includes Conjecture \ref{Euler}, has been
proposed, and proved in some special cases, in \cite[Conj. 1.4]{Sokal}.
We feel
that the following stronger property (which could be called ``monotone total 
positivity'') actually holds. 
\begin{con} \label{EulerMTP}
For all $i,j,r \in \NN$  
the determinant of the submatrix determined by the rows indexed by 
$i,i+1, \ldots , i+r$ and columns indexed by $j,j+1, \ldots , j+r$
is a monotonically increasing function of $i \in \NN$.	
\end{con}
\noindent 
We have checked that this is true if $j+r \leq 44$ and $i+r \leq 44$.
Note that Conjecture \ref{EulerMTP} implies Conjecture \ref{Euler}
by \cite[Thm. 2.8]{Pin}.

\medskip 
{\bf Acknowledgments:} 
Some of the computations for the research presented in this paper have been carried out using some Maple packages for computing with Coxeter systems and 
posets developed by Pietro Mongelli and John Stembridge.
I would like to thank Mario Marietti and Volkmar Welker for pointing out some relevant references.
The author was partially supported by the MIUR Excellence Department Projects 
CUP E83C18000100006 and E83C23000330006.

\bibliographystyle{amsplain}

\end{document}